\documentclass[12pt,letterpaper]{article}
\usepackage[utf8]{inputenc}
\usepackage[T1]{fontenc}
\usepackage{lmodern}
\usepackage[margin=1.0in]{geometry}
\usepackage{microtype}
\usepackage[english=usenglishmax]{hyphsubst}
\usepackage[numbers,super]{natbib}
\usepackage{amsmath,mathtools,amsfonts,amssymb}
\usepackage[pdftex]{graphicx}
\usepackage{epstopdf}
\usepackage[hyphens]{url}
\usepackage[linesnumbered,ruled]{algorithm2e}
\usepackage{booktabs}
\usepackage{colortbl}
\usepackage{xcolor}
\usepackage{hyperref}
\usepackage{setspace}      %singlespacing/doublespacing[onehalfspacing]
\usepackage{blindtext}

\title{On the performance of preconditioned methods to solve \(L^p\)-norm phase unwrapping}
\author{Ricardo Legarda-Saenz, Carlos Brito-Loeza, Arturo Espinosa-Romero\\
CLIR at Facultad de Matemáticas, Universidad Autónoma de Yucatán\\
Apartado Postal 172. 97110 Mérida, Yucatán. México\\
E-mail: \texttt{rlegarda@correo.uady.mx}  } 
\date{\today}

\begin{document}
\maketitle 

\begin{abstract}
In this paper, we analyze and evaluate suitable preconditioning techniques to improve the performance of the $L^p$-norm phase unwrapping method. We consider five preconditioning techniques commonly found in the literature, and analyze their performance with different sizes of wrapped-phase maps.
\end{abstract}

\section{Introduction}
There has been a growing interest in the development of new techniques for the processing of coherent signals. This kind of signals is generated by measurement techniques like synthetic aperture radar (SAR), magnetic resonance imaging (MRI), and interferometry among others. The objective of this processing is to estimate the phase term $\phi_{\mathbf{x}}$, also known as phase map, from a signal which general model could be expressed as \[U_{\mathbf{x}} = S_{\mathbf{x}}\exp\left( \mathit{i}\phi_{\mathbf{x}}\right),\]  where $\mathbf{x} = (x,y)$ is the position and $S_{\mathbf{x}}$ is the signal amplitude. The estimation of the phase term $\phi_{\mathbf{x}}$ becomes very relevant, given that this term can be related to different physical quantities such as geographical topography in the case of SAR, or the optical path difference in the case of optical interferometry~\cite{Gasvik2002,Tupin2014}.

However, as it can be seen from the signal model, the estimation of the phase term $\phi_{\mathbf{x}}$ is not straightforward. Instead, the estimated phase map from the signal is defined as 
\begin{equation}\label{eq:wrapped}
\psi_{\mathbf{x}} = \phi_{\mathbf{x}} + 2\pi k_{\mathbf{x}},
\end{equation}
where $k_{\mathbf{x}}$ is a function that bounds the values to $-\pi < \psi_{\mathbf{x}}\leq\pi.$ The term $\psi_{\mathbf{x}}$ is known as wrapped phase and is a nonlinear function of $\phi_{\mathbf{x}}.$ This term is not useful for measurements because just offers the principal values of the phase term $\phi_{\mathbf{x}},$ so it is necessary estimate $\phi_{\mathbf{x}}$ from this wrapped phase $\psi_{\mathbf{x}}$. This process is called phase unwrapping~\cite{Ghiglia1998}.

Phase unwrapping is an ill-posed problem~\cite{Bertero1998,Vogel2002}. The unwrapping process consists of integrating the gradient field of the wrapped phase map~\cite{Ghiglia1998}. Even in ideal conditions, phase estimation is not trivial due to the non linearity of $\psi_{\mathbf{x}}$ caused by cyclic discontinuities. In real conditions, the unwrapping process becomes very difficult, where noise signal, sub-sampling or differences larger than $2\pi$ (real or not) generate ambiguities hard to process not allowing the accurate recovering of the phase map $\phi_{\mathbf{x}}.$ 

In the literature there are two main strategies to solve the unwrapping problem~\cite{Ghiglia1998,Gens2003}: the first one, known as path-following or local methods, consists of integrating the differences of the wrapped phase over a path that covers the entire phase map. The second one, considers the problem globally and the solution is presented in the form of an integral over the wrapped region. In the global strategy, there are two approaches: the first consists of using Green functions together with a numerical solution based on the Fourier transform~\cite{Fornaro1996a,Fornaro1996,Lyuboshenko1999,Lyuboshenko2000}. 

The second approach expresses the solution as a $L^p$-norm minimization problem, resulting on the solution of weighted differential equations. The numerical solution of these differential equations leads to a nonlinear system $\mathbf{A}(\mathbf{u})\mathbf{u} = \mathbf{b},$ which has to be solved iteratively with great computational cost~\cite{Ghiglia1996,Hooper2007,Guo2014a}. Typically, conjugate gradient (CG) or multigrid method are the methods of choice for this kind of problems. In the case of the conjugate gradient, matrix $\mathbf{A}$ is expected to be well conditioned, otherwise the convergence will be slow. However, very often one comes across with ill conditioned matrices when working with these nonlinear systems. The solution is to precondition the matrix $\mathbf{A}$; this is, instead of solving the original system $\mathbf{A}(\mathbf{u})\mathbf{u} = \mathbf{b},$ we solve the preconditioned system $\mathbf{M}^{-1}\mathbf{A}(\mathbf{u})\mathbf{u} = \mathbf{M}^{-1}\mathbf{b}.$  The matrix $\mathbf{M}$, called a preconditioner for the matrix $\mathbf{A}$, is chosen to improve the condition number of the matrix $\mathbf{A}.$ In most cases, this preconditioning matrix is problem dependent. 

The goal of this paper is to analyze and evaluate suitable preconditioning techniques to improve the performance of the $L^p$-norm phase unwrapping method~\cite{Ghiglia1996,Ghiglia1998}. We consider five preconditioning techniques commonly found in the literature, and analyze their performance with different sizes of wrapped-phase maps. The organization of this paper is as follows: first, we describe the $L^p$-norm phase unwrapping method and their numerical solution. Then, the performance of the selected preconditioning techniques is evaluated by numerical experiments with different sizes of a synthetic wrapped-phase map. Finally, we discuss our results and present some concluding remarks.

\section{Methodology}
The $L^p$-norm based method proposed by D. C. Ghiglia and M.D. Pritt for 2-D phase unwrapping is defined as~\cite{Ghiglia1996}
\begin{equation}\label{eq:funcional}
\underset{\phi}{\min}\; J(\phi_{\mathbf{x}}) = \iint_{\Omega}\left\lvert\frac{\partial\phi}{\partial x} - \frac{\partial\psi}{\partial x}\right\rvert^{p} + \left\lvert\frac{\partial\phi}{\partial y} - \frac{\partial\psi}{\partial y}\right\rvert^{p}\;d\mathbf{x}
\end{equation}
 where $\mathbf{x} = (x,y),$ and $\Omega \subseteq \mathbb{R}^{2}$ is the domain of integration.
To obtain the solution of the problem expressed in Eq.  (\ref{eq:funcional}), the first-order optimality condition or Euler-Lagrange equation has to be derived, resulting in the following partial differential equation (PDE)
\begin{equation}\label{eq:gradiente}
- \frac{\partial}{\partial x}\left[\left(\frac{\partial\phi}{\partial x} - \frac{\partial\psi}{\partial x}\right)\left\lvert\frac{\partial\phi}{\partial x} - \frac{\partial\psi}{\partial x}\right\rvert^{p-2} \right] 
- \frac{\partial}{\partial y}\left[\left(\frac{\partial\phi}{\partial y} - \frac{\partial\psi}{\partial y}\right)\left\lvert\frac{\partial\phi}{\partial y} - \frac{\partial\psi}{\partial y}\right\rvert^{p-2}\right] = 0,
\end{equation}
with boundary conditions
\[\left[ \left(\frac{\partial\phi}{\partial x} - \frac{\partial\psi}{\partial x}\right)\left\lvert\frac{\partial\phi}{\partial x} - \frac{\partial\psi}{\partial x}\right\rvert^{p-2}, \; \left(\frac{\partial\phi}{\partial y} - \frac{\partial\psi}{\partial y}\right)\left\lvert\frac{\partial\phi}{\partial y} - \frac{\partial\psi}{\partial y}\right\rvert^{p-2} \right] \cdot \hat{\mathbf{n}} = 0,\] where $\hat{\mathbf{n}}$ denotes the unit outer normal vector to the boundary. 

\subsection{Numerical solution}
Let $u_{i,j} = u(x_i,y_j)$ to denote the value of a function $u_{\mathbf{x}}$ at point $(x_i,y_j)$ defined on $\Omega = [a,b]\times[c,d],$ where the sampling points are $x_i = a + (i-1)h_x,\quad y_j = c + (j-1)h_y,$ with $1 \leq i \leq M,\; 1 \leq j \leq N,$ $h_x = (b-a)/(M-1),\; h_y = (d-c)/(N-1)$ and $M, N$ are the number of points in the discrete grid of points. We use $u$ to represent any of the variables $\phi$ and $\psi$ defined in the previous equations. Derivatives are approximated using standard forward and backward finite difference schemes \[\delta_{x}^{\pm}u_{i,j} = \pm\frac{u_{i,j\pm 1}-u_{i,j}}{h_x}\quad\text{and}\quad\delta_{y}^{\pm}u_{i,j} = \pm\frac{u_{i\pm 1j}-u_{i,j}}{h_y}.\]  The gradient and the divergence are approximated as \[\triangledown  u_{i,j} = (\delta_{x}^{+}u_{i,j},\delta_{y}^{+}u_{i,j})\quad\text{and}\quad\triangledown\cdot\triangledown  u_{i,j} = \delta_{x}^{-} (\delta_{x}^{+}u_{i,j}) + \delta_{y}^{-}(\delta_{y}^{+}u_{i,j}),\] respectively. 

Hence the numerical approximation of the Eq. (\ref{eq:gradiente}) is given by 
\begin{equation}\label{eq:gradientep}
- \delta_{x}^{-}\left[\left(\delta_{x}^{+}\phi - \Delta_{i,j}^{x}\right)\left\lvert\delta_{x}^{+}\phi - \Delta_{i,j}^{x}\right\rvert^{p-2}\right] - \delta_{y}^{-}\left[\left(\delta_{y}^{+}\phi - \Delta_{i,j}^{y}\right)\left\lvert\delta_{y}^{+}\phi - \Delta_{i,j}^{y}\right\rvert^{p-2}\right] = 0,
\end{equation}
with boundary conditions 
\[\left[ \left(\delta_{x}^{+}\phi - \Delta_{i,j}^{x}\right)\left\lvert\delta_{x}^{+}\phi - \Delta_{i,j}^{x}\right\rvert^{p-2}, \; \left(\delta_{y}^{+}\phi - \Delta_{i,j}^{y}\right)\left\lvert\delta_{y}^{+}\phi - \Delta_{i,j}^{y}\right\rvert^{p-2} \right] \cdot \hat{\mathbf{n}} = 0,\]
where the terms $\Delta_{i,j}^{x}$ y $\Delta_{i,j}^{y}$ are defined as
\[\Delta_{i,j}^{x} = 
  \begin{cases}
    \mathcal{W}\left(\psi_{i,j+1} - \psi_{i,j}\right)\:&\text{if}\quad 1 \leq i \leq M-1,\:1 \leq j \leq N-2\\
    \mathcal{W}\left(\psi_{i,j} - \psi_{i,j-1}\right)\:&\text{if}\quad 1 \leq i \leq M-1,\: j = N-1\\
  \end{cases}
\]
\[\Delta_{i,j}^{y} = 
  \begin{cases}
    \mathcal{W}\left(\psi_{i+1,j} - \psi_{i,j}\right)\:&\text{if}\quad 1 \leq i \leq M-2,\:1 \leq j \leq N-1\\
    \mathcal{W}\left(\psi_{i,j} - \psi_{i-1,j}\right)\:&\text{if}\quad i = M-1,\:1 \leq j \leq N-1\\
  \end{cases}
\]
where $\mathcal{W}$ is the wrapping operator~\cite{Ghiglia1998}. Applying the previous discrete approximations, we have that the numerical solution of Eq. (\ref{eq:gradiente}) is given by
\begin{equation}\label{eq:ecuacion}
  \begin{aligned}
  &\left(\phi_{i,j} - \phi_{i,j-1} - \Delta_{i,j-1}^{x}\right)U_{i,j-1} - \left(\phi_{i,j+1} - \phi_{i,j} - \Delta_{i,j}^{x}\right)U_{i,j}\\
  + &\left(\phi_{i,j} - \phi_{i-1,j} - \Delta_{i-1,j}^{y}\right)V_{i-1,j} - \left(\phi_{i+1,j} - \phi_{i,j}  - \Delta_{i,j}^{y}\right)V_{i,j}
  = 0,
  \end{aligned}
\end{equation}
with boundary conditions 
\[\;\left[\left(\phi_{i,j+1} - \phi_{i,j} - \Delta_{i,j}^{x}\right)U_{i,j},\;\left(\phi_{i+1,j} - \phi_{i,j}  - \Delta_{i,j}^{y}\right)V_{i,j}\right]\cdot \hat{\mathbf{n}} = 0,\] 
where for simplicity and without loss of generality we consider $h_x = h_y  = 1.$ The terms $U$ y $V$ are defined as
\begin{equation}\label{eq:UV}
  \begin{aligned}
  U_{i,j} &= \left\lvert\phi_{i,j+1} - \phi_{i,j} - \Delta_{i,j}^{x}\right\rvert^{p-2} = \frac{\tau }{\left\lvert\phi_{i,j+1} - \phi_{i,j} - \Delta_{i,j}^{x}\right\rvert^{2-p} + \tau},\\
  V_{i,j} &= \left\lvert\phi_{i+1,j} - \phi_{i,j} - \Delta_{i,j}^{y}\right\rvert^{p-2} = \frac{\tau}{\left\lvert\phi_{i+1,j} - \phi_{i,j} - \Delta_{i,j}^{y}\right\rvert^{2-p} + \tau}.
  \end{aligned}
\end{equation}
to force the values of the terms to lie in the range (0,1). This helps the stability and convergence of the numerical solution~\cite{Bloomfield1983,Ghiglia1996}. Usually, $\tau = 0.01$ is selected for phase unwrapping process~\cite{Ghiglia1996,Ghiglia1998}. 

Discretization of Eq. (\ref{eq:gradiente}) leads to a nonlinear PDE because the terms $U$ and $V$ are functions of the input data and the solution. This is solved using the following iterative procedure~\cite{Bloomfield1983,Scales1988,Ghiglia1996}: first, given an initial value $\phi,$ the terms $U$ and $V$ are computed; then, the terms $U$ and $V$ are held fixed and Eq. (\ref{eq:ecuacion}) is solved using preconditioned conjugate gradient~\cite{Shewchuk1994,Golub1996}. With the current solution $\phi,$ the terms $U$ and $V$ are updated and a new solution is computed. This process is repeated until convergence. 

Now, we present the algorithm to solve the discretization of Eq. (\ref{eq:gradiente}). First we arrange Eq. (\ref{eq:ecuacion}) in  matrix form $\mathbf{A}\phi = \mathbf{b}$ as
\begin{equation}\label{eq:Axb}
  \begin{aligned}
  &-\left(\phi_{i,j+1}U_{i,j} + \phi_{i+1,j}V_{i,j} + \phi_{i,j-1}U_{i,j-1} + \phi_{i-1,j}V_{i-1,j}\right)\\
  &+\left(U_{i,j-1} + U_{i,j} + V_{i-1,j} + V_{i,j}\right)\phi_{i,j}\\
  &= \Delta_{i,j-1}^{x}U_{i,j-1} - \Delta_{i,j}^{x}U_{i,j} + \Delta_{i-1,j}^{y}V_{i-1,j} - \Delta_{i,j}^{y}V_{i,j}.
  \end{aligned}
\end{equation}
Notice that $\mathbf{A}$ is a $MN\times MN$ sparse matrix which depends on the terms $U$ and $V,$ so it needs to be constructed at each iteration. The preconditioned conjugate gradient (PCG) used in our work is the implementation proposed in Ref. \cite{Shewchuk1994}, and the explicit structure  of the algorithm is given in Algorithm \ref{algoritmo}.
\begin{algorithm}
  \DontPrintSemicolon
  \KwData{the wrapped phase $\psi$ and $p < 2.$}
  \KwResult{the unwrapped phase $\phi.$}
  \BlankLine
  $k \leftarrow 0,\;error \leftarrow 1$\;
  $\phi_{i,j}^{k}\leftarrow\text{random values}$\;
  \While{$k < k_{max}\;\mathbf{ and}\; error > tol $}
  {
    compute $U$ and $V$, Eq. (\ref{eq:UV})\;
    solve Eq. (\ref{eq:Axb}) using PCG:\;
    \Begin
    {
      construct $\mathbf{A}$ and $\mathbf{b}$ \;
      estimate preconditioning matrix $\mathbf{M}$ from $\mathbf{A}$\;\label{matrizM}
      $\phi^{k+1} \leftarrow\phi^{k}$\;
      $l \leftarrow 0$\;
      $\mathbf{r}\leftarrow\mathbf{b} - \mathbf{A}\phi^{k+1}$\;
      $\mathbf{d}\leftarrow\mathbf{M}^{-1}\mathbf{r}$\;
      $\delta_{new}\leftarrow\mathbf{r}^{T}\mathbf{d}$\;
      $\delta_{0}\leftarrow\delta_{new}$\;
      \While{ $l < l_{max}\;\mathbf{ and}\; \delta_{new} > \epsilon^{2}\delta_{0}$}
        {
          $\mathbf{q}\leftarrow\mathbf{A}\mathbf{d}$\;
          $\alpha\leftarrow\delta_{new} / \mathbf{d}^{T}\mathbf{q}$\;
          $\phi^{k+1} \leftarrow\phi^{k+1} + \alpha\mathbf{d}$\;
          \eIf{ $l$ is divisible by $\sqrt{MN}$ }
            {$\mathbf{r}\leftarrow\mathbf{b} - \mathbf{A}\phi^{k+1}$\;}
            {$\mathbf{r}\leftarrow\mathbf{r} - \alpha\mathbf{q}$\;} 
          $\mathbf{s}\leftarrow\mathbf{M}^{-1}\mathbf{r}$\;
          $\delta_{old}\leftarrow\delta_{new}$\;
          $\delta_{new}\leftarrow\mathbf{r}^{T}\mathbf{s}$\;
          $\beta\leftarrow\delta_{new} / \delta_{old}$\;
          $\mathbf{d}\leftarrow\mathbf{s} + \beta\mathbf{d}$\;
          $l\leftarrow l +1$\;  
        }
    }
    $error\leftarrow\|\phi_{i,j}^{k+1}-\phi_{i,j}^{k}\| / \|\phi_{i,j}^{k}\|$\;    
    $k\leftarrow k +1$\;  
  }
  \caption{\(L^p\)-norm phase unwrapping algorithm.}\label{algoritmo}
\end{algorithm}

The key point of the performance of the algorithm is the proper selection of the preconditioning matrix $\mathbf{M}$, step \ref{matrizM} of Algorithm \ref{algoritmo}. $\mathbf{M}$ can be defined in many different ways but it should meet the following requirements: a) the preconditioned system should be easy to solve, and b) the preconditioning matrix should be computationally cheap to construct and apply~\cite{Saad2003}. 

Here we test suitable preconditioning techniques to estimate matrix $\mathbf{M}$ in Algorithm \ref{algoritmo} and analyze the performance of the $L^p$-norm phase unwrapping method. For this purpose, we selected the following commonly found preconditioning techniques in literature~\cite{VanderVorst2003,Saad2003,Chen2005a}:
\begin{enumerate}
  \item $\mathbf{M}\Leftarrow\mathbf{I},$ where $\mathbf{I}$ is the identity matrix; with this matrix, PCG becomes CG~\cite{Golub1996}.
  \item $\mathbf{M}\Leftarrow\mathbf{D},$ where $\mathbf{D}$ is the diagonal of matrix $\mathbf{A};$ this technique is known as Jacobi preconditioning.
  \item $\mathbf{M}\Leftarrow ILU(\mathbf{A}),$ the incomplete LU factorization with no fill-in.
  \item $\mathbf{M}\Leftarrow IC(\mathbf{A}),$ the incomplete Cholesky factorization with no fill-in.
  \item $\mathbf{M}\Leftarrow SOR(\mathbf{A}),$ the successive over-relaxation factorization.
\end{enumerate}
The main reason for selecting these preconditioning techniques was their no fill-in property. That is, the preconditioning matrices preserve the sparse structure of matrix $\mathbf{A}$ and require just the same memory space used by matrix $\mathbf{A}.$

\section{Numerical experiments}
To illustrate the performance of the selected preconditioning techniques, we carried out the numerical experiments using a Intel\textsuperscript{\tiny\textregistered} Core\textsuperscript{\tiny\texttrademark} i7 @ 2.40 GHz laptop with Debian GNU/Linux 10 (buster) 64-bit and 16 GB of memory. For our experiments, we programmed all the functions using C language, GNU g++ 8.3 compiler and Intel\textsuperscript{\tiny\textregistered} MKL 2019 library. It is important to highlight that all the functions were programmed from scratch and were programmed from the basic algorithms, without using modifications that improve their performance.

\begin {figure}[!ht]
  \begin{center}
    \includegraphics[width=0.7\textwidth,keepaspectratio]{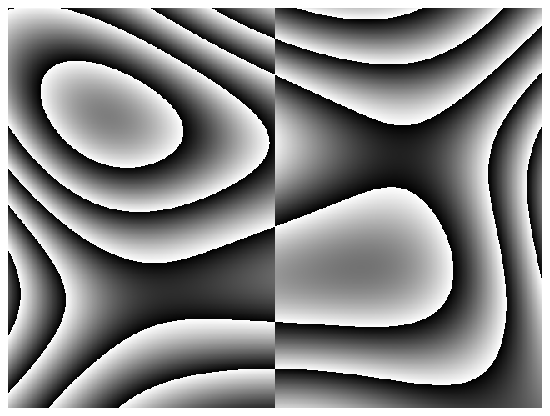}
  \end{center}
  \caption{Wrapped phase map used in the experiments.}\label{fig:Wrapped}
\end {figure}

We use the wrapped phase map shown in Fig. \ref{fig:Wrapped} as the  data term $\psi$ of the Algorithm \ref{algoritmo}. This wrapped phase map has a resolution of 640x480 pixels, and this resolution was used as our reference scale.

We generated several scaled wrapped phase map from the one shown in Fig. \ref{fig:Wrapped}, and used them as data in Algorithm \ref{algoritmo}. The scaled image sizes and the resultant sizes of sparse matrix $\mathbf{A}$ in our experiments are shown in Table \ref{tab:table}.
\begin{table}[!ht]
  \caption{Image sizes used in the experiments.}\label{tab:table}
  \begin{center}
  \renewcommand{\arraystretch}{1.5}
  \setlength{\tabcolsep}{6pt}
    \begin{tabular}{ccrc}\toprule
      &&\multicolumn{2}{c}{matrix $\mathbf{A}$}\\ \cline{3-4}
      scale&image size&variables&density (\%)\\ 
      \midrule
      \rowcolor{black!20}0.25& 160 x 120&95440&0.0250\\
      0.50& 320 x 240&382880&0.0064\\
      \rowcolor{black!20}0.75&480 x 360& 862320&0.0028\\
      1.00&640 x 480&1533760&0.0016\\
      \rowcolor{black!20}1.25&800 x 600&2397200&0.0010\\
      1.50&960 x 720&3452640&0.0007\\
      \rowcolor{black!20}1.75&1120 x 840&4700080&0.0005\\
      2.00&1280 x 960&6139520&0.0004\\
     \bottomrule
    \end{tabular}
  \end{center}
\end{table}

\begin {figure}[!ht]
  \begin{center}
    \includegraphics[width=\textwidth,keepaspectratio]{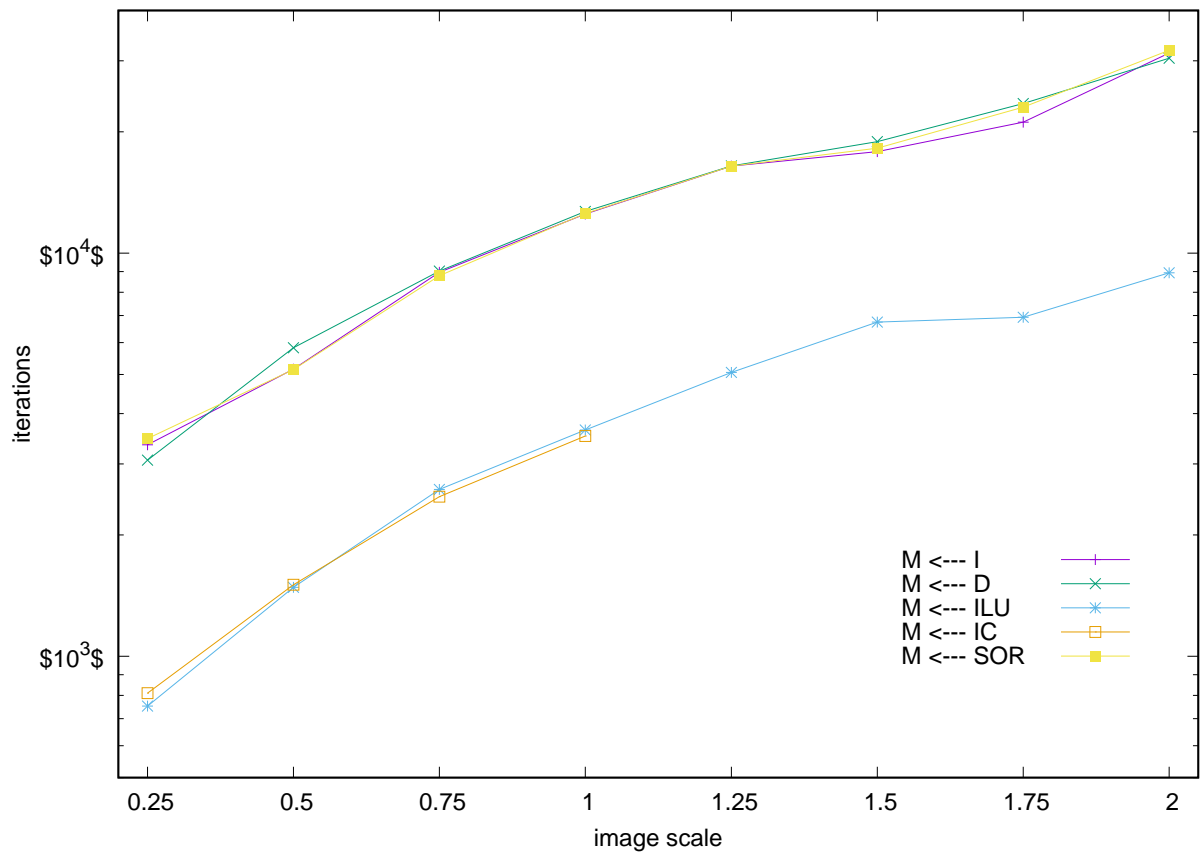}
  \end{center}
  \caption{Iterations employed for the different test.}\label{fig:Itera}
\end {figure}

\begin {figure}[!ht]
  \begin{center}
    \includegraphics[width=\textwidth,keepaspectratio]{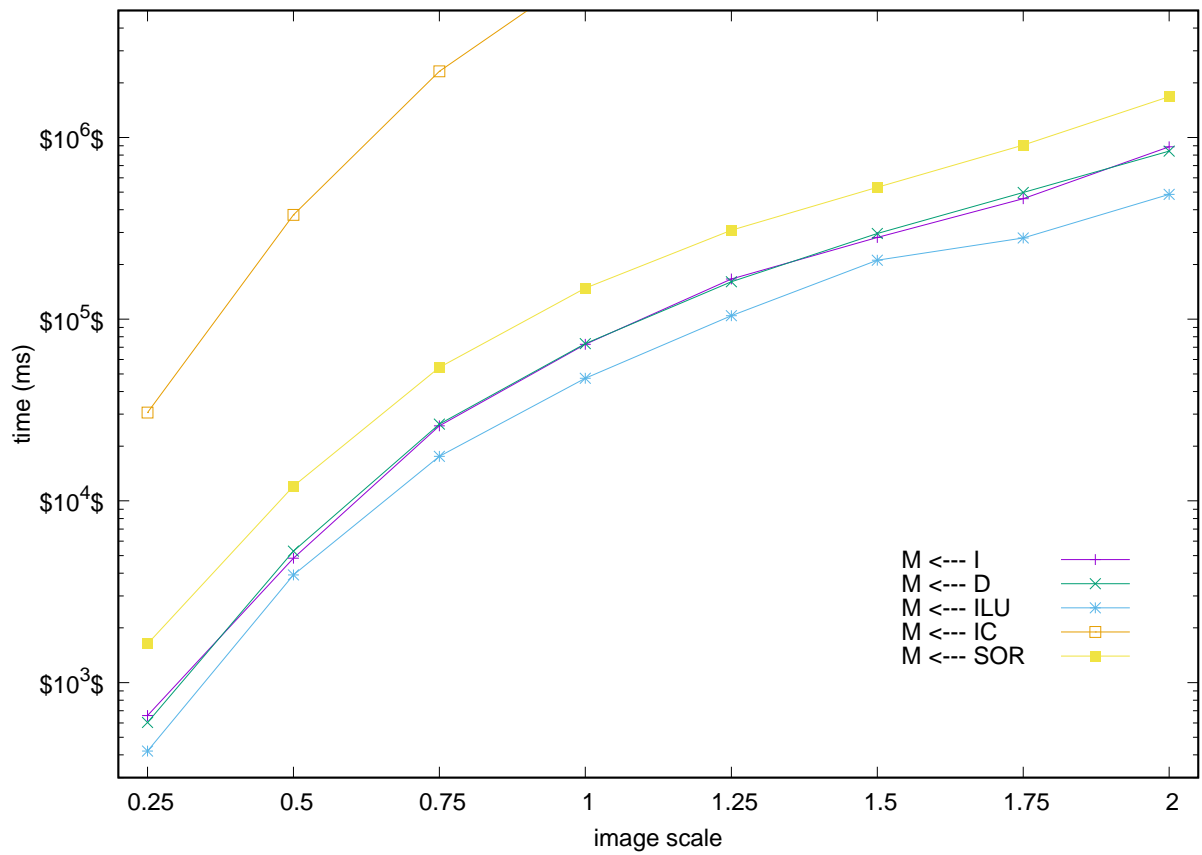}
  \end{center}
  \caption{Computational time used for the different test.}\label{fig:Tiempo}
\end {figure}

For each scaled wrapped phase map,  we tested the selected five preconditioning matrices.  The stopping criteria  used in Algorithm \ref{algoritmo} were $k_{max} = 500,$ $tol = 10^{-6},$ $l_{max} = 2MN,$ $\epsilon = 0.005.$ Fig. \ref{fig:Itera} shows the iterations needed to obtain the solution for the different scaled wrapped phase maps, and Fig. \ref{fig:Tiempo} shows the computational time consumed for each preconditioning technique.
It is worth the value to remark that  the time taken to construct preconditioning matrix $\mathbf{M}$: except for the incomplete Cholesky factorization, all the preconditioning techniques employed between 1\% and 2\% of the total processing time. However, the incomplete Cholesky factorization took more than 90\%. This is why in Figs 2 and 3, we only show the first four experiments with this technique.

\begin {figure}[!ht]
  \begin{center}
    \includegraphics[width=0.7\textwidth,keepaspectratio]{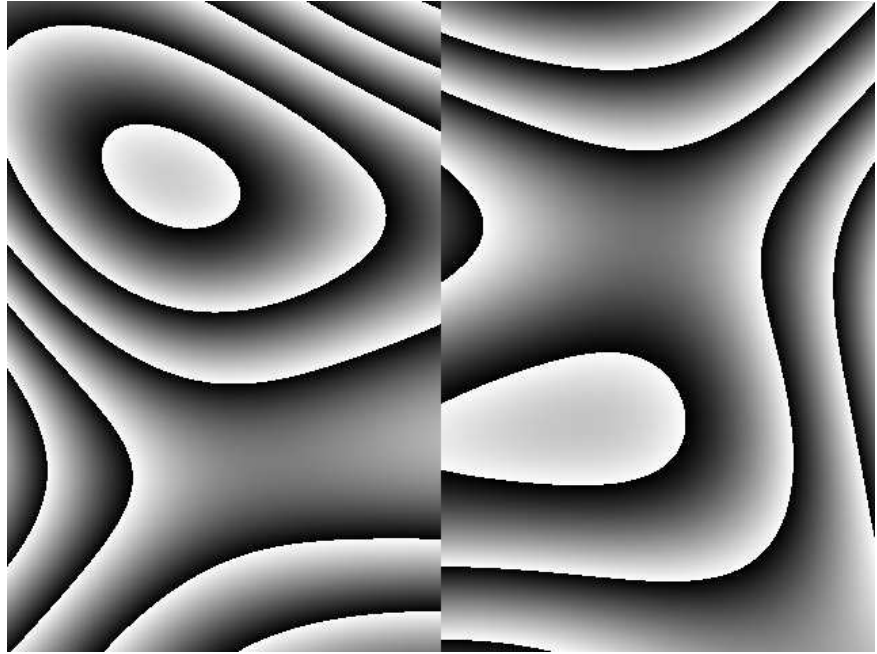}
  \end{center}
  \caption{Resultant unwrapped phase map using Algorithm \ref{algoritmo}. The phase maps was wrapped for purposes of illustration.}\label{fig:Unwrapped}
\end {figure}

Finally, we use a normalized error $Q$ to compare the unwrapping estimation; this error is defined as~\cite{Perlin2016}:
\begin{equation}
Q\left(\mu,\nu\right) = \frac{\| \mu - \nu\|_{2}}{\| \mu\|_{2} + \|\nu\|_{2}}, 
\end{equation}
where $\mu$ and $\nu$ are the signals to be compared. The normalized error values vary between zero (for perfect agreement) and one (for perfect disagreement). For all the cases, we found that the normalized error was around $Q = 0.17$. An example of the resultant unwrapped phase map is shown in Fig. \ref{fig:Unwrapped}.

\section{Discussion of results and conclusions}
From the results obtained, we have the following remarks. First, any of the techniques used had no impact on the obtained results, since for all the experiments the normalized error Q was approximately the same, $Q = 0.17$. The differences are evident when we analyze the performance in terms of iterations and computational time. If we only consider the number of iterations, the incomplete LU factorization  and incomplete Cholesky factorization clearly show their advantage over the other three methods. However, this perception changes when we also consider the computational time used. Clearly, the incomplete Cholesky factorization consumed a lot of time which makes it unviable to be considered as a preconditioning technique. In general, numerical results show the incomplete LU factorization as the best choice to be used as preconditioning technique in Algorithm 1. 

As work in the future, we are going to analyze the use of dedicated libraries for estimating the preconditioning matrix, where we hope to obtain better results.

%%%%%%%%%%%%%%%%%%%%%%%%%%%%%%%%%%%%%%%%%%%%%%%%%%%
%%%    Referencias
%%%%%%%%%%%%%%%%%%%%%%%%%%%%%%%%%%%%%%%%%%%%%%%%%%%
\bibliographystyle{splncs04}
%\bibliography{/home/rlegarda/Dropbox/Articulos/MendeleyDropbox/bibtex/library}

%%%%%%%%%%%%%%%%%%%%%%%%%%%%%%%%
%           Figuras
%%%%%%%%%%%%%%%%%%%%%%%%%%%%%%%%
\newpage
\section*{Figures}

\end{document}